\documentclass[12pt]{article}
\usepackage{amsmath}
\usepackage{amsfonts}
\usepackage{amssymb}
\usepackage{theorem}

\title{Generic groups and the weak amalgamation property}
\author{A. Ivanov and K. Majcher }

\newtheorem{theorem}{Theorem}[section]
\newtheorem{proposition}[theorem]{Proposition}
\newtheorem{corollary}[theorem]{Corollary}
\newtheorem{lemma}[theorem]{Lemma}
\newtheorem{definition}[theorem]{Definition}
\newtheorem{example}[theorem]{Example}
\newtheorem{remark}[theorem]{Remark}

\begin{document} 
\topmargin = 12pt
\textheight = 630pt 
\footskip = 39pt 

\maketitle

\begin{quote}
{\bf Abstract} 
We consider the logic space of countable (enumerated) groups and show that closed subspaces corresponding to some standard classes of groups have (do not have) generic groups. 
We also discuss the cases of semigroups and associative rings. 
\end{quote}

\section{Introduction} 
The recent paper \cite{GEL} considers countable groups as elements of a Polish space in the following way. 
An {\em enumerated group} is the set $\omega$ together with a multiplication function 
$\cdot : \omega \times \omega \rightarrow \omega$, 
an inversion function  
$^{-1}: \omega \rightarrow \omega$ , and an identity element $e \in  \omega$ defining a group. 
The space $\mathcal{G}$ of enumerated groups is a closed subset of the space 
$\mathcal{X}=\omega^ { \omega \times \omega} \times \omega^{\omega} \times \omega$ (under the product topology). 
Any abstract property of groups, say $P$, naturally defines the invariant subset $\mathcal{G}_P$ of $\mathcal{G}$. 
Assuming that $\mathcal{G}_P$ is Polish, it is studied in \cite{GEL} which group-theoretic properties are generic in $\mathcal{G}_P$.  
In particular, are there generic groups 
in $\mathcal{G}_P$? 
(A group $G$ is {\em generic} if all isomorphic copies of $G$ form a comeagre subset of $\mathcal{G}_P$.) 
Among the rich variety of results of that paper we mention as the most notable the statement that there is no generic group in $\mathcal{G}$, and  moreover, when $P$ is the property of left orderability, the space $\mathcal{G}_P$ does not have a generic group. 
The authors deduce this using involved arguments in the style of model-theoretic forcing.   
\footnote{Some of these results were independently (but later) obtained in \cite{marton}. 
This approach was also applied in \cite{marton-top}.}

A. Ivanov studied the problem of existence of generic objects in \cite{iva99} in full generality. 
He introduced some general condition which is now called the {\em weak amalgamation property (WAP)}
and showed that it together with the joint embedding property characterizes the existence of generics. 
\footnote{ 
A. Ivanov used in \cite{iva99} the name the almost amalgamation property. 
This condition was considered later in the case of automorphisms in the very influential paper \cite{KR} where the name WAP was used. } 

The goal of our paper is as follows. 
Using Ivanov's approach we build a universal framework for investigations of generic objects in various classes of enumerated groups, semigroups and rings. 
Our meta-theorems are direct and short. 
In particular, see Theorem \ref{nWAP} which describes a very general property implying the absence of generic objects.  
The framework is applied in a number of cases. 
Firstly, we deduce the results of \cite{GEL} mentioned above in a short fashion. 
Furthermore, we prove: 
\begin{itemize} 
\item there is no generic group in any Burnside variety for odd exponents greater than $C$ (where $C$ is some natural number), see Corollary \ref{burnside};
\item  there is no generic semi-group, see Corollary \ref{darbisemi};  
\item there is no generic associative ring, see Corollary \ref{ring}. 
\end{itemize} 
We also apply our methods in the contrary direction, i.e. for proving of existence of generics in some natural group varieties. 
In particular, in Section 2.3 we prove the following theorem. 
\bigskip

\noindent 
{\bf Theorem 2.9.} 
{\em Let $c,p\in \omega \setminus \{0,1\}$ and $p$ be prime $>\mathsf{max}(2,c)$. 
Then the space of nilpotent groups of degree $c$ and of exponent $p$ is a closed subspace of $\mathcal{G}$ which  
has a generic group. } 
  
\bigskip 

\noindent   
The framework also works in the case of actions of groups on first-order structures. 
In particular, we show that there is no generic group action on $(\mathbb{Q},<)$, see Corollary \ref{Q}. 
This result is a contribution to the topic of generic representations of countable groups, see \cite{DM}. 

\subsection{General preliminaries}
Fix a countable atomic structure $M$ in a language $L$. 
Let $T$ be an expansion of $Th(M)$ in some $L\cup \bar{r}$ where $\bar{r}= (r_1,\ldots,r_{\iota})$ is a sequence of additional relational/functional symbols. 
We assume that $T$ is axiomatizable by sentences of the following form: 
$$
(\forall \bar{x}) \left( \bigvee_{i} (\phi_i(\bar{x}) \land \psi_i(\bar{x}) ) \right) , 
$$
where $\phi_i (\bar{x})$ is a quantifier-free formula of the language $L \cup \bar{r}$, and $\psi_i (\bar{x})$ is a first order formula of the language $L$ (in particular, $Th(M)\subseteq T$). 
Let $\mathcal{X}_{M,T}$ be the space of all $\bar{r}$-expansions of $M$ to models of $T$. 
This is a topological space with respect to so called {\em logic topology}. 
In the following few paragraphs we introduce some notions previously defined in \cite{iva99}. 
They will also be helpful in a description of the logic topology. 
\begin{remark} 
{\em In \cite{iva99} it is assumed that $Th(M)$ is $\omega$-categorical. 
It is worth mentioning that the arguments used in the proof of Theorem 1.2 in \cite{iva99} work under the weaker assumption that $M$ is atomic. 
However, in the present paper we will work under the assumption of $\omega$-categoricity. 
}
\end{remark}  

For a finite $A \subset M$ we define a {\em diagram} of $\bar{r}$ on $A$ as follows. 
To every functional symbol $r_i$ from $\bar{r}$, say of arity $\ell$,  a partial function 
$\mathsf{r}_i :A^{\ell} \to A$ is associated. 
When $r_i$ is a relational symbol, we include into the diagram a single formula from every pair 
$\{ \mathsf{r}_i(\bar{a}'), \lnot \mathsf{r}_i(\bar{a}') \}$, 
where $\bar{a}'$ is a tuple from 
$A$ of the corresponding length. 
We will say that $A$ is the
{\em domain} of this diagram. 
For ease of notation, we will usually view $A$ as a tuple of pairwise distinct elements, say $\bar{a}$. 
Then given a diagram $D$ on $A$, we write 
$\bar{a} = \mathsf{Dom}(D)$.  
It will also be convenient to use the form $D(\bar{a})$, where the domain is already given in the expression. 
    
Given a diagram $D(\bar{a})$, $\bar{a} \subset M$, 
which is satisfied in some $(M, {\bf \bar{r}}) \models T$, we view $D(\bar{a})$ as the  atomic diagram of the corresponding partial substructure of $(M, {\bf \bar{r}})$. 

Let $\mathcal{B}_T$ be the set of all theories 
of the form $D(\bar{a})\cup T \cup Th(M,\bar{a})$ 
that are satisfied in some expansions $(M, {\bf \bar{r}}) \models T$. 
Since $D(\bar{a})$ uniquely defines the corresponding member of $\mathcal{B}_T$, 
we will identify this member with the diagram. 
Furthermore, since $M$ is atomic, each element of $\mathcal{B}_T$ is determined by a formula of the form $\phi (\bar{a}) \wedge \psi (\bar{a})$, where $\psi (\bar{x})$ is a complete formula for $M$ realized by $\bar{a}$, and 
$\phi (\bar{a})$ is a quantifier-free formula in the language $L \cup \bar{r}$. 
The corresponding $\phi(\bar{x}) \wedge \psi (\bar{x})$ is called {\em basic}. 

For every diagram  $D(\bar{a}) \in \mathcal{B}_T$ the set 
\[ 
[D] = \{ (M, {\bf \bar{r}}) \in \mathcal{X}_{M,T} \, | \, (M,{\bf \bar{r}}) 
\mbox{ satisfies } D(\bar{a})\} 
\] 
is clopen in the logic topology. 
In fact, the family 
$\{ [D] \, | \, D (\bar{a})\in \mathcal{B}_T \}$ is usually taken as a basis of this topology. 
It is metrizable. 
Indeed, take an enumeration of $M^{<\omega}$, say 
$M^{<\omega}= \{ \bar{a}_0 ,\bar{a}_1 ,\ldots \}$. 
Define a metric $d$ metrizing $\mathcal{X}_{M,T}$ as follows: 
\[ 
d((M,{\bf \bar{r}}) , (M,{\bf \bar{s})} ) = \sum \{ 2^{-n} \, : \, {\bf \bar{r}}, \, {\bf \bar{s}} \mbox{ do not agree on } \bar{a}_n  \} .
\] 
Let $D\in \mathcal{B}_T$ and $\mathsf{Dom}(D) = \bar{a}$. 
Then $D$ can be viewed as an expansion of $M$ by finite relations corresponding to $\bar{r}$. 
When $\mathsf{r}_i$ is a functional symbol, the relation is $\mathsf{Graph}(\mathsf{r}_i )$, the graph of the corresponding partial function on $\bar{a}$. 
For any automorphism 
$\alpha \in \mathsf{Aut}(M)$ we define  $\alpha (D)$ to be the diagram obtained from $D(\bar{a})$ by substitution $\alpha(\bar{a})$ instead of $\bar{a}$.    

The set $\mathcal{B}_T$ is ordered by 
extension, it will be denoted by $\subseteq$.  
More formally: 
$D \subseteq D'$ if $\mathsf{Dom}(D) \subseteq \mathsf{Dom}(D')$ 
and $D'$ implies $D$ under $T$ (in particular, the partial functions defined in $D'$ extend the corresponding partial functions defined in $D$). 
Note that when $\bar{a}$ is a tuple  of some entries of $\bar{b}$, $D'(\bar{b})$ implies $D(\bar{a})$ under $T$ exactly when so do the corresponding basic formulas. 
 
Since the space $\mathcal{X}_{M,T}$ is Polish, the Baire category method can be applied. 
\begin{definition} 
An expansion $(M, {\bf \bar{r}})\models T$ is called {\em generic} if it has a comeagre isomorphism class in $\mathcal{X}_{M,T}$. 
\end{definition} 
In order to present the description of  generic expansions from \cite{iva99} we need definitions of JEP, AP, CAP and WAP. 
\begin{itemize} 
\item The family $\mathcal{B}_T$ has the {\em joint embedding property} if for any two elements 
$D_1, D_2 \in \mathcal{B}_T$ there is $D_3\in\mathcal{B}_T$ 
and an automorphism $\alpha \in \mathsf{Aut}(M)$ such that $D_1 \subseteq D_3$ and $\alpha (D_2) \subseteq D_3$. 
\item The family $\mathcal{B}_T$ has the {\em amalgamation property} if for any $D_0 ,D_1 $, $D_2 \in \mathcal{B}_T$ with $D_0 \subseteq D_1$ and $D_0 \subseteq D_2$ there is $D_3 \in \mathcal{B}_T$ and an automorphism $\alpha \in \mathsf{Aut}(M)$ fixing $\mathsf{Dom} (D_0 )$ such that $D_1 \subseteq D_3$ and $\alpha (D_2) \subseteq D_3$. 
\item The family $\mathcal{B}_T$ has the {\em cofinal amalgamation property} if for any $D_0 \in \mathcal{B}_T$ there is an extension $D'_0 \in \mathcal{B}_T$ such that for any $D_1 , D_2 \in \mathcal{B}_T$ with $D'_0 \subseteq D_1$ and $D'_0 \subseteq D_2$ there is $D_3 \in \mathcal{B}_T$ and an automorphism $\alpha \in \mathsf{Aut}(M)$ fixing $\mathsf{Dom} (D'_0)$ such that $D_1 \subseteq D_3$ and $\alpha (D_2) \subseteq D_3$. 
\item The family $\mathcal{B}_T$ has the {\em weak amalgamation property} if for any $D_0 \in \mathcal{B}_T$ there is an extension $D'_0 \in \mathcal{B}_T$ such that for any $D_1 , D_2 \in \mathcal{B}_T$ with $D'_0 \subseteq D_1$ and $D'_0 \subseteq D_2$ there is $D_3 \in \mathcal{B}_T$ and an automorphism $\alpha \in \mathsf{Aut}(M)$ fixing $\mathsf{Dom} (D_0)$ such that $D_1 \subseteq D_3$ and $\alpha (D_2) \subseteq D_3$. 
\end{itemize} 
By Theorem 1.2 from \cite{iva99}  
\begin{quote} 
$\mathcal{X}_{M,T}$ {\em has a generic expansion $(M,{\bf \bar{r}})$ if and only if the family $\mathcal{B}_{T}$ has JEP and WAP. }
\end{quote}

\subsection{The case of groups} 

\paragraph{I. The basic case} 
The structure $M$ is just $\omega$ and $L$ consists of one constant symbol $1$ which is interpreted by number $1$.  
Let $\bar{r}$ consist of the binary function of multiplication $\cdot$ and a unary function $^{-1}$. 
Furthermore, a finite sequence of constants can be  included into $\bar{r}$; this does not change the exposition. 
Let $T$ be the universal theory of groups with the unit $1$   
\footnote{It is not necessary to fix $1$ in the language of $M$. 
For example, the paper \cite{GEL} does not do this. 
We fix it just for convenience of notations in Scenario {\em II} below.}. 
 
Each element of $\mathcal{B}_T$ consists of a tuple $\bar{c}\subset \omega$ containing $1$ and partial functions for $\cdot$ and $^{-1}$ on $\bar{c}$.  
{\em This is a finite partial group. } 
(When $\bar{r}$ has additional constants, it is assumed that distinguished constants are included into $\bar{c}$.) 
The space $\mathcal{X}_{M,T}$ is the logic space of all groups of size $\aleph_0$, where the unit is fixed. 
{\em We emphasize that this is exactly the space $\mathcal{G}$ of enumerated groups from \cite{GEL}.}  
Indeed, Corollary 3.1.3 of \cite{GEL} gives a basis of clopen sets of $\mathcal{G}$ as follows. 
Take a finite system of group equations and inequations $\Sigma (\bar{x})$ and a tuple $\bar{a} \subset \omega$ with $|\bar{x}| = |\bar{a}|$. 
Let $[ \Sigma (\bar{a}) ]$ be the set of all group operations on $\omega$ for which the tuple $\bar{a}$ satisfies $\Sigma (\bar{x})$. 
Then the family of all sets of the form 
$[\Sigma (\bar{a})]$ is the basis of $\mathcal{G}$. 
On the other hand, note that every diagram $D(\bar{c}) \in \mathcal{B}_T$ determines a finite system, say $\Sigma_D (\bar{c})$, of group equations and inequations over $\bar{c}$ (for words which values are defined under the corresponding restrictions of the group operations) so that $[D]=[\Sigma_D(\bar{c})]$.   
Furthermore, each $[\Sigma (\bar{a})]$ is a union of basic clopen sets determined by some family of diagrams $D(\bar{c}) \in \mathcal{B}_T$. 
Indeed, consider all finite partial multiplication tables on tuples $\bar{c} \subset \omega$ with $\bar{a}\subseteq \bar{c}$ such that every subword $u(\bar{a}')$ of a word from $\Sigma (\bar{a})$ has a value in $\bar{c}$ and $\Sigma (\bar{a})$ 
is satisfied under this table. 
Each multiplication table of this form is viewed as a diagram $D(\bar{c})\in \mathcal{B}_T$. 

\bigskip 

If instead of $T$ above one takes a universal extension of the theory of groups, say $U$, then the corresponding space of enumerated groups is a closed subspace of $\mathcal{G}$. 
Such a subspace will be denoted by $\mathcal{G}_{U}$. 
 
In the paragraph below we will have several sorts in $M$. 
The sort just described will always occur. 
It will be denoted by $\mathsf{Gp}$. 
The symbols $\bar{r}$ will always contain ones  for multiplication and for inversion. 
In particular, if $(M, {\bf \bar{r}})\in \mathcal{X}_{M,T}$ then we denote by $\mathsf{Gp_{\bar{r}}}$ the group defined by the corresponding operations from ${\bf \bar{r}}$ on $\mathsf{Gp}$.

\paragraph{II. The case of an action} 
Let $M_0$ be a countable atomic structure of some language $L_0$. 
Define $M$ to be $M_0$ with an additional sort $\omega$ called $\mathsf{Gp}$ and the constant symbol $1$ interpreted by $1$. 
The symbols $\bar{r}$ include $\cdot$,  $^{-1}$ and a new symbol $\mathsf{ac}$ for a function $M_0\times \mathsf{Gp}\rightarrow M_0$. 
Now the theory $T$ contains the universal axioms of groups on the sort $\mathsf{Gp}$ with the unit $1$.  
We also add the following axioms for an action:  
\[ 
\mathsf{ac}(x, z_1 \cdot z_2 ) = \mathsf{ac}(\mathsf{ac}(x,z_1 ),z_2 ) \, \, , \, \,  \mathsf{ac}(x, 1 )= x. 
\] 
The space $\mathcal{X}_{M,T}$ is the logic space of all countable expansions of $M$ where the group structure is defined on $\mathsf{Gp}$, with a fixed unit $1$ and an action $\mathsf{ac}$ . 

It is natural to add the universal axioms that $\mathsf{ac}$ preserves the structure of $M_0$. 
In this way we obtain the space of actions on $M_0$ by automorphisms.

\subsection{Generics are generic} 

In this section we indicate a connection of our approach (and approaches of \cite{GEL}, \cite{marton} and \cite{marton-top}) with model-theoretic forcing. 
Let $\mathcal{X}$ be a Polish space with a countable basis of clopen sets $\mathcal{A}$.  
For $A\in \mathcal{A}$ and a Borel set $B$ we write 
$A \Vdash B$ if $A\setminus B$ is meagre in $A$, and $A\Vdash \sim B$ if $A\cap B$ is meagre in $A$. 
We say that $A$ {\em decides} $B$ if exactly one of these cases holds.  
If $\mathcal{B}$ is a family of Borel subsets, then following Definition 2.2 of \cite{MI}, we say that an element $x\in \mathcal{X}$ is $\mathcal{B}$-{\em generic} in $\mathcal{X}$, if  for every $B\in \mathcal{B}$ there is a neighborhood $A\in \mathcal{A}$ of $x$ that decides $B$. 
 
This notation is strongly connected with model-theoretic forcing. 
Recall that given a countable language $L$, a set of constant symbols $C$ and an $L$-theory $T$, a {\em condition} is a finite set $p$ of $(L\cup C)$-formulas without free variables which are atomic or negations of atomic formulas, such that $p$ is consistent with $T$.     
We will assume that $C=\omega$ and $T$ is universal. 
It is clear that each condition $p$ 
defines a set $[ p ]$ of all $L$-structures defined on $\omega$ which satisfy $p \cup T$. 
This is a clopen set in the logic space $\mathcal{X}_{T}$ of all models of $T$ defined on $\omega$.   
Let $\mathcal{A}$ be the basis of clopen set of $\mathcal{X}_T$ consisting of such $[p ]$. 

We now follow Section 8 of \cite{keisler}.  
When $\phi$ is a formula over $L\cup C$ without free variables, then  Robinson's forcing $p \Vdash \phi$ is defined as follows:   

- $\phi \in p$, when $\phi$ is atomic,

- there is no extension $q\supseteq p$ with $q \Vdash \psi$, when $\phi = \neg \psi$, 

- there is $c\in \omega$ with $p\Vdash \psi (c)$, when $\phi = \exists x \psi$, 

- $p \Vdash \psi_1$ or $p\Vdash \psi_2$, when $\phi = \psi_1 \vee \psi_2$. \\ 
The following statement is folklore (and can be easily verified): 
\begin{quote} 
Let $B= [\phi ]$ be the Borel set of all $L$-structures on $\omega$ which satisfy $\phi\cup T$. 
Then $p\Vdash \phi$ implies $[p] \Vdash B$.   
\end{quote} 
We now define generic structures as follows. 
Let $M$ be a structure on $\omega$ and 
$D(M)$ be its diagram of atomic sentences and their negations. 
The structure $M$ is called {\em generic} for $T$ if every finite subset of $D(M)$ is a condition for $T$ and for every sentence $\phi$ of $L \cup C$ there a finite subset $p\subseteq D(M)$ such that $p \Vdash \phi$ or $p\Vdash \neg \phi$. 
It is well known that every $T$-generic structure is existentially closed. 

Now consider the logic space $\mathcal{G}_U$ for some universal theory $U$. 
This is a Polish $Sym (\omega)$-space under the action by permutations on $\omega$. 
According to the previous paragraph if $M$ is a generic structure for $U$ then for every Borel set of the form $[\phi ]$ there is a basic open $A$ of the form $[ p ]$ with $M \in [p ]$ such that $[p ]$ decides $[\phi ]$. 
As we already know, this property of $M$ is exactly $\mathcal{B}$-{\em genericity} in the terminology of \cite{MI}, where in our case $\mathcal{B}$ is the family of all Borel subsets of $\mathcal{G}_U$ of the form $[\phi ]$ for first-order $\phi$ 
\footnote{the paper \cite{MI} also considers other examples of families $\mathcal{B}$}. 
By Proposition 2.4 of \cite{MI} all $\mathcal{B}$-generic elements of $\mathcal{G}_U$ form a dense $G_{\delta}$-subset of $\mathcal{G}_U$. 
This proves the following statement. 

\begin{proposition} 
The set of $U$-generic groups is an invariant comeagre subset of $\mathcal{G}_U$. 
The same property holds for the set of existentially closed groups. 
\end{proposition} 
Since algebraic closedness coincides with existential closedness, this proposition strengthens Lemma 5.2.7 of \cite{GEL}, Theorem 6.1 of 
\cite{marton} and Theorem 3.11 of \cite{marton-top}. 
Summarizing this section, we conclude: {\em if $G$ is a generic group of the space of enumerated groups $\mathcal{G}_U$, then $G$ is $U$-generic and existentially closed.}

\section{Seeking generics} 

In this section, we give examples of group varieties $U$ such that the corresponding $\mathcal{G}_U$ contains generic groups. 
In particular, we prove Theorem \ref{cnilp} formulated in the introduction. 
In Section 2.1 we introduce the crucial notion of {\em t-isolation} and show how it works for proving WAP.  

\subsection{t-Isolation}  
\noindent 
Consider a countable atomic structure $M$ in a language $L$. 
Let $T$ be an expansion of $Th(M)$ in some $L\cup \bar{r}$ where $\bar{r}= (r_1,\ldots,r_{\iota})$ is a sequence of additional relational/functional symbols. 
We assume that $T$ satisfies all the assumptions stated in the beginning of Section 1.1. 
Let $M$ have a distinguished sort $\mathsf{Gp}$. 
Under Scenarios $I$ or $II$ of Section 1.2, it is assumed that $\mathsf{Gp}$ is a sort for a group.  
However, in Definition \ref{isol} we do not take this assumption. 

Let $\bar{c}$ be a tuple of the sort $\mathsf{Gp}$ of the structure $M$. 
For any expansion $(M,{\bf \bar{r}})$ the tuple $\bar{c}$ determines a theory, say $\Xi_{\, \bar{c}, {\bf \bar{r}}}$,  which consists of all valid conditions $r_i (t_1 (\bar{c}_1), \ldots , t_k(\bar{c}_k ))$ (or their negations) and $t_1 (\bar{c}') = (\not= ) t_2(\bar{c}'')$, where $\bar{c}_1, \ldots , \bar{c}_k$, $\bar{c}'$ and $\bar{c}''$ are subtuples of $\bar{c}$ and $t_1 (\bar{x}_1), \ldots, t_k(\bar{x}_k )$, $t_1 (\bar{x}')$ and $t_2 (\bar{x}'')$ are terms of $L\cup \bar{r}$. 
In the following definition we distinguish the case when $\Xi_{\, \bar{c}, \bf \bar{r}}$ is uniquely determined by some diagram. 

\begin{definition} \label{isol}
Let $\bar{c}$ be a subtuple of $\bar{d}$ and $D(\bar{d}) \in \mathcal{B}_T$ (and we do not assume that $\bar{c}$ is of the sort $\mathsf{Gp}$). 
We say that $D(\bar{d})$ is {\em t-isolating} (term-isolating) {\em for} $\bar{c}$ if 
the theory 
$\Xi_{\, \bar{c}\cap \mathsf{Gp}, {\bf \bar{r}}}$ is the same in every member $(M,{\bf \bar{r}}) \in [D(\bar{d})]$. 
\end{definition} 

\noindent 
In this definition $\bar{c}\cap \mathsf{Gp}$
denotes the subtuple of $\bar{c}$ of all entries which belong to $\mathsf{Gp}$. 

When $(M,{\bf \bar{r}})$ is an expansion of $M$ as in the definition, we may consider $\mathsf{Gp}$ as a structure with respect to those relations/functions of $L\cup \bar{r}$ which do not depend on elements of other sorts. 
Then any subset of $\mathsf{Gp}$ closed under 
operations of $\mathsf{Gp}$ of this kind is viewed as a substructure. 
It is worth noting here, that the theory $\Xi_{\, \bar{c}, \bf \bar{r}}$ defined before the definition is the same with the corresponding theory defined in this substructure.  

Definition \ref{isol} will be especially important in cases {\em I - II} of Section 1.2.  
In that cases, we may assume that the theory $\Xi_{\, \bar{c}\cap\mathsf{Gp}, {\bf \bar{r}}}$ consists of all valid conditions $w(\bar{c}') = (\not= ) 1$ where $\bar{c}'$ is a subtuple of $\bar{c}\cap \mathsf{Gp}$ and $w(\bar{x}')$ is a group word. 
Then it is worth noting that Definition \ref{isol} implies that the tuple $\bar{c}$ generates the same marked group $\langle \bar{c} \cap \mathsf{Gp} \rangle$ in all members $(M, {\bf \bar{r}})\in [D(\bar{d})]$  (as a subgroup of the group defined on $\mathsf{Gp}$).  

From now on in this section we assume the circumstances of the basic case,  i.e. case {\em I}. 
In particular, the structure $M$ has only one sort $\mathsf{Gp}$. 
Then, 
$\Xi_{\, \bar{c}\cap \mathsf{Gp}, {\bf \bar{r}}}= \Xi_{\, \bar{c}, {\bf \bar{r}}}$ in the definition.

Let $G$ be a subgroup of an enumerated group and $G = \langle \bar{c} \rangle$ 
for some tuple $\bar{c} \subseteq \omega$. 
We say that $G$ is {\em t-isolated} by $D(\bar{d})$, if  $D(\bar{d})$ is t-isolating for $\bar{c}$, and $G$ is the corresponding subgroup $\langle \bar{c} \rangle$. 

\begin{example} \label{2_2}
{\em Let $F= \langle \bar{c} \rangle$ be a finite group which extends to an enumerated group of $T$. 
Let $\bar{d}$ be an enumeration of $F$. 
The multiplication table of $F$ together with inequalities between distinct elements, can be viewed as a diagram $D(\bar{d})$ which is t-isolating for $\bar{c}$. 
} 
\end{example} 

\noindent 
Now it is easy to see that every finite group which can be embedded into an enumerated group of $T$ is a t-isolated group in the sense of the following definition.   

\begin{definition}  
We say that an abstract group $G$ is {\em t-isolated with respect to } $T$ if there is an enumerated group, say $G_1$, which is defined on $\mathsf{Gp}$ (under Case {\em I}) and there is a tuple $\bar{c} \subset \omega$ such that the subgroup $\langle \bar{c} \rangle \le G_1$,  is isolated by some diagram $D(\bar{c}_1)$ for $\bar{c}$, where $\bar{c}_1 \subseteq \langle \bar{c} \rangle$ in every member of $[D(\bar{c}_1)]$. 
\end{definition} 

\noindent 
Note that when the condition above holds, then extending the diagram if necessary, it can be arranged that $\bar{c}_1$ is presented by the same tuple of words on $\bar{c}$ in all members of $[D(\bar{c}_1)]$, and, furthermore,  $\bar{c}$ is a subtuple of $\bar{c}_1$.  
Another easy observation is the fact that the set of all enumerated groups from $\mathcal{X}_{M,T}$ which contain a t-isolated $G$ (with a fixed set of generators $\bar{c}$), form a clopen subset of $\mathcal{X}_{M,T}$. 

\begin{definition}
We say that t-isolated diagrams are {\em dense} in $\mathcal{B}_T$ if any $D_0 (\bar{c}) \in \mathcal{B}_T$ extends to some $D(\bar{c}_1)$ which is t-isolating for $\bar{c}$ and 
$\bar{c}_1 \subseteq \langle \bar{c} \rangle$ in every member of $[D(\bar{c}_1)]$. 
\end{definition}

The following proposition will be our tool for seeking generics. 

\begin{proposition} \label{WAPgrp}
Under the circumstances of case I assume that t-isolated diagrams are dense in $\mathcal{B}_T$. 

Then WAP for $\mathcal{B}_T$ is equivalent to the following property: 
\begin{quote} 
every t-isolated $G_0$ can be extended to a t-isolated $G_1$ such that any two t-isolated extensions of $G_1$ can be amalgamated over $G_0$. 
\end{quote} 
\end{proposition} 

{\em Proof.} 
Necessity. 
Assuming that $G_0$ has a copy, 
say $\langle \bar{c}_0 \rangle$, which is t-isolated by some $D_0 (\bar{c}'_0)$ with  
$\bar{c}'_0 \subseteq \langle \bar{c}_0 \rangle$,   
find an extension $D(\bar{c})\supseteq D_0 (\bar{c}'_0)$ as in the formulation of WAP. 
Extend it to some t-isolating $D'(\bar{c}')$ for $\bar{c}$ with $\bar{c}' \subseteq \langle \bar{c} \rangle$. 
It is easy to see that $D' (\bar{c}')$ witnesses WAP for $D_0(\bar{c}'_0)$. 
Let $G_1$ be the corresponding group $\langle \bar{c}\rangle$.  
Having $G_2 =\langle \bar{c}_2 \rangle$ and $G_3 =\langle \bar{c}_3 \rangle$ which contain $G_1$ and are t-isolated by some $D'_2 (\bar{c}'_2)$ and $D'_3 (\bar{c}'_3)$ with $\bar{c}'_2 \subseteq \langle \bar{c}_2 \rangle$ and $\bar{c}'_3 \subseteq \langle \bar{c}_3 \rangle$ (and where we may assume $\bar{c} \subseteq \bar{c}_2 \cap \bar{c}_3$), 
we apply WAP for amalgamation $D'_2 (\bar{c}'_2)$ and $D'_3 (\bar{c}'_3)$ over $D_0 (\bar{c}'_0)$. 
Assume that $D_4 (\bar{c}_4)$ is the corresponding amalgam. 
Then $\langle \bar{c}_4  \rangle$ is a group containing $\bar{c}_2 \cup \bar{c}_3$.  
By t-isolation, $\bar{c}_2$ and $\bar{c}_3$ generate subgroups of $\langle \bar{c}_4 \rangle$ which are copies of $G_2$ and $G_3$. 

Sufficiency. 
Let $D (\bar{c})\in \mathcal{B}_T$. 
Extend it to a t-isolating $D_0 (\bar{c}')$ for $\bar{c}$ with $\bar{c}' \subseteq \langle \bar{c} \rangle$. 
Let $G_0 = \langle \bar{c} \rangle$ be the corresponding t-isolated group.  
Find $G_1 = \langle \bar{c}_1 \rangle$ as in the formulation and the corresponding t-isolating diagram $D_1 (\bar{c}'_1)$ with 
$\bar{c}'_1 \subseteq \langle \bar{c}_1 \rangle$. 
We may assume that $\bar{c}'\subseteq\bar{c}_1$.  
In order to verify amalgamation of extensions of $D_1 (\bar{c}'_1)$ we first note that by density of t-isolating diagrams, we may restrict ourselves only to extensions that are t-isolating for some subtuples of their domains.  
If $D_2 (\bar{c}'_2)$ and $D_3 (\bar{c}'_3)$
are such extensions and $G_2 = \langle \bar{c}_2 \rangle$ and $G_3 = \langle \bar{c}_3 \rangle$  are the corresponding groups with $\bar{c}'_2 \subseteq \langle \bar{c}_2 \rangle$ and $\bar{c}'_3 \subseteq \langle \bar{c}_3 \rangle$, then by t-isolation, $G_1$ is a subgroup of both $G_2$ and $G_3$. 
Find $G_4$ which amalgamates $G_2$ and $G_3$ over $G_0$. 
It obviously belongs to some $[D_4 (\bar{c}_4)]$ where  $D_4 (\bar{c}_4)$ amalgamates  $D_2 (\bar{c}_2)$ and $D_3 (\bar{c}_3)$ over $D_0 (\bar{c}')$. 
$\Box$

\begin{remark} 
{\em It is worth noting that a slightly weaker form of the density assumption follows from WAP. 
Indeed, assume that there is $D(\bar{c}) \in \mathcal{B}_T$ such that for every $D'(\bar{c}') \in \mathcal{B}_T$ extending $D(\bar{c})$ the diagram $D'(\bar{c}')$ is not t-isolating for $\bar{c}$ (where $\bar{c}' \subseteq \langle \bar{c} \rangle$ is not assumed). 
Therefore, for such $D'(\bar{c}')$ two distinct groups can be realized as $\langle \bar{c} \rangle$ in some members of $[D'(\bar{c}')]$. 
Thus $\mathcal{B}_T$ does not have WAP.  
}
\end{remark}

The following remark will be very useful below.

\begin{remark} \label{RF} 
{\em 
As we already know, every finite group which extends to an enumerated group of $T$, is t-isolated (see Example \ref{2_2}). 
Using this we can conclude that  {\em if $T$ defines a variety of groups and all finitely generated groups satisfying $T$ are residually finite, then t-isolated diagrams are dense in $\mathcal{B}_T$}. 
Indeed, assuming that $D(\bar{c})\in \mathcal{B}_T$ take a finite group $F$ generated by $\bar{c}$ and satisfying all relations which appear in $D(\bar{c})$.  
Let $\bar{c}_1$ be an enumeration of $F$. 
Then the multiplication table of $F$ (together with inequalities between distinct elements) is viewed as the corresponding diagram $D(\bar{c}_1)$, where  we also have $\bar{c}_1 \subseteq \langle \bar{c} \rangle$. 
The diagram $D(\bar{c}_1)$ works as an appropriate t-isolating extension of $D(\bar{c})$. 
}
\end{remark} 

We also add here that under the assumptions of Remark \ref{RF},  t-isolated groups satisfying $T$ are exactly finite ones. 
Indeed, if $G= \langle \bar{c} \rangle$ is infinite, then any basic set $[D'(\bar{c}')]$ containing it, also contains a finite group. 

\subsection{Abelian groups} 

Abelian groups form a closed subspace of $\mathcal{G}$. 
We denote it by $\mathcal{G}_{Abel}$ 
where $Abel$ denotes the corresponding axiomatization. 
Under the scenario of case {\em I} of Section 1.2 we have 
$\mathcal{G}_{Abel}= \mathcal{X}_{M,Abel}$ . 

The following theorem is proved in \cite{marton} and \cite{marton-top}. 
We give an easy and independent proof of it. 

\begin{theorem} 
The space $\mathcal{G}_{Abel}$ has a generic group. 
\end{theorem} 

{\em Proof.} 
According to Theorem 1.2 of \cite{iva99} we have to verify that JEP and WAP hold for $\mathcal{B}_{Abel}$.  
Since abelian groups form a variety (we call it $\mathbb{G}_{Abel}$), JEP is satisfied: just apply $\oplus$ for domains of the corresponding diagrams. 
In order to verify WAP we will apply Proposition \ref{WAPgrp}. 

Firstly, note that since finitely generated abelian groups are residually finite, Remark \ref{RF} implies that the t-isolated diagrams are dense in $\mathcal{B}_{Abel}$. 

Let us verify the property that every t-isolated $G_0$ can be extended to a t-isolated $G_1$ such that any two t-isolated extensions of $G_1$ can be amalgamated over $G_0$.  
As we already mentioned after Remark \ref{RF}, t-isolated abelian groups are exactly finite ones. 
Now the key point is the fact that the variety $\mathbb{G}_{Abel}$ has the amalgamation property. 
Indeed, if $A_0$ is a subgroup of $A_1$ and $\phi$ is an embedding of $A_0$ into $A_2$, then the group $(A_1 \oplus A_2)/\{ (x -\phi (x)) | x\in A_0 \}$ is the corresponding amalgamation. 
$\Box$

\subsection{Nilpotent groups of prime odd exponent}

Let $c,p\in \omega \setminus \{0,1\}$ and $p$ be prime $>\mathsf{max}(2,c)$. 
Nilpotent groups of degree $c$ and of exponent $p$ form a closed subspace of $\mathcal{G}$. 
We denote it by $\mathcal{G}_{c_{nil}p}$. 

\begin{theorem} \label{cnilp}
The space $\mathcal{G}_{c_{nil}p}$ has a generic group. 
\end{theorem} 

{\em Proof.} 
Consider $\mathcal{G}_{c_{nil}p}$ as the space $\mathcal{X}_{M,c_{nil}p}$ 
under the scenario of case {\em I} of Section 1.2.   
According to Theorem 1.2 of \cite{iva99} we have to verify that JEP and WAP hold for $\mathcal{B}_{c_{nil}p}$.  
Since nilpotent groups of degree $c$ and of exponent $p$ form a variety (we call it $\mathbb{G}_{c_{nil}p}$), JEP is satisfied: apply $\oplus$ for domains of the corresponding diagrams. 
In order to verify WAP we will apply Proposition \ref{WAPgrp}. 

First, note that since finitely generated nilpotent groups are residually finite (see Section 5.2 in \cite{CMZ}), the t-isolated diagrams are dense in $\mathcal{B}_{c_{nil}p}$. 
Furthermore, t-isolated groups of the variety $\mathbb{G}_{c_{nil}p}$ are exactly finite ones, see Remark \ref{RF} and the comment after it. 
Thus, Proposition \ref{WAPgrp} can be applied to finite groups. 
However, we need some additional algebraic material. 

Let $G \in \mathbb{G}_{c_{nil}p}$ and 
\[ 
\langle 1 \rangle = \Gamma_{c+1}(G) \leq \Gamma_{c}(G) \le  \ldots \le \Gamma_{1}(G) =G ,
\]
\[ 
\langle 1 \rangle = Z_{0}(G) \leq Z_{1}(G) \leq \ldots \leq Z_{c}(G)=G 
\]
be the lower and the upper central series of $G$ respectively.
It is clear that 
\begin{quote}
for any group extension $A \leq B$ and for any $i\le c$ we have  
$\Gamma_{i}(A)\le \Gamma_{i}(B) \cap A $. 
\end{quote}
We now need the following lemma. 

\begin{lemma} \label{lcnilp} 
For any finite group 
$G_0 \in \mathbb{G}_{c_{nil}p}$ there is a finite group $G_1 \in \mathbb{G}_{c_{nil}p}$ such that $G_0 \le G_1$ and for every finite 
$C \in \mathbb{G}_{c_{nil}p}$ containing $G_1$ and every $i \le c$, we have that 
$\Gamma_{i}(G_1) \cap G_0 = \Gamma_{i}(C)\cap G_0$. 
\end{lemma}

{\em Proof.} 
Given $k\le c$, we show how to build $G_1$ that satisfies the conclusion of the lemma for all $i\le k$. 
To satisfy the conclusion of the lemma for $k=1$, take $G_1 = G_0$. 
We now apply induction. 
Let $k<c$ and  assume that we already have some $H \ge G_0$ satisfying the statement for all $i \le k$. 
Let $G_1$ be a finite extension of $H$ with maximal $\Gamma_{k+1} (G_1 ) \cap G_0$.  
Since for all $i\le k$ we have 
$\Gamma_{i}(G_1) \cap G_0 = \Gamma_{i}(H)\cap G_0$,   
the statement of the conclusion of the lemma holds for $G_1$ and all $i\le k+1$. 
$\Box$ 

\bigskip 

In \cite{baud04} A. Baudisch defines expansions of groups from $\mathbb{G}_{c_{nil}p}$ by additional unary predicates $P_2, \ldots, P_c$. 
In order to present these expansions let 
$G$ be a group in $\mathbb{G}_{c_{nil}p}$ and let $Q_1 (G) = P_1 (G) = G$. 
When $1 < n \le c + 1$ let 
$Q_n (G) = \langle \bigcup_{\ell +k= n} [P_{\ell} (G),P_k (G)] \rangle$ and $P_n (G)$  be a subgroup of $G$ with 
$Q_n (G) \subseteq P_n (G) \subseteq Z_{c+1-n} (G)$ 
and $P_n (G) \subseteq  P_{n-1}(G)$. 
Let $\mathbb{G}^P_{c,p}$ be the class of all groups in $\mathbb{G}_{c_{nil}p}$ with additional predicates $P_2(G), \ldots, P_c(G)$ as described above. 
Let $L^+$ be the corresponding language that extends the language of groups. 
Note that if we extend a group $G \in \mathbb{G}_{c_{nil}p}$ by predicates  
$Q_i (G)=P_i (G) =\Gamma_i (G)$ for 
$1 < i \le c$ then we get an $L^+$-group from $\mathbb{G}^P_{c,p}$.
The following statement is the main result of 
\cite{baud04} 
(see Theorem 1.2 and the discussion after it). 
\begin{quote}
The class of finite structures from $\mathbb{G}^P_{c,p}$ has the hereditary property and the amalgamation property.
\end{quote}

We can now finish the proof of Theorem \ref{cnilp}. 
Let us verify the amalgamation condition from 
Proposition \ref{WAPgrp}. 
Let $G_0 \in \mathbb{G}_{c_{nil}p}$ be finite, and let $G_1 \in \mathbb{G}_{c_{nil}p}$ be a finite group as in Lemma \ref{lcnilp}. 
Let $G_2, G_3 \in \mathbb{G}_{c_{nil}p}$ be two finite extensions of $G_1$.
 
Consider $G_1$, $G_2$ and $G_3$ as $L^+$-structures defining $P_i$ on them by $\Gamma_i (G_1 )$, $\Gamma_i (G_2 )$ and $\Gamma_i (G_3 )$ respectively ($1 < i \le c$).    
Now for any $i \leq c$ let 
$P_i(G_0)=\Gamma_i(G_1) \cap G_0$. 
By the hereditary property $(G_0, P_2(G_0), \ldots, P_c(G_0))$ is also an $L^+$-structure from $\mathbb{G}^P_{c,p}$.

Since $G_1$ is a subgroup of $G_2$ and $G_3$ then by Lemma \ref{lcnilp}, for any $i\le c$, we have that 
$P_i(G_0)=P_i(G_1) \cap G_0 = P_i(G_2) \cap G_0=  P_i(G_3) \cap G_0$. 
In particular $(G_0 , P_2 (G_0 ), \ldots , P_c (G_0)  )$ is an $L^+$-substructure of both 
$(G_2 , P_2 (G_2 ), \ldots , P_c (G_2)  )$ and $(G_3 , P_2 (G_3 ), \ldots , P_c (G_3)  )$. 
By Theorem 1.2 of \cite{baud04} there is $L^+$-structure $G_4$ amalgamating $G_2$ and $G_3$ over $G_0$. 
Thus, the condition of Proposition \ref{WAPgrp} is satisfied.
$\Box$ 

\bigskip

We finish this section by a discussion if WAP in Proposition \ref{WAPgrp} can be replaced by CAP, the cofinal amalgamation property. 
We recommend \cite{KKKP} for some general material concerning this question. 
Note that in the case $c=2$ the following statement strengthens the final argument of the proof of Theorem \ref{cnilp}. 

\begin{proposition} \label{nil-cap}
The variety $\mathbb{G}_{2_{nil}p}$ has the following property:  
any finite $G_0$ can be extended to a finite $G_1$ such that every two finite extensions of $G_1$ can be amalgamated over $G_1$ in a finite group.  
\end{proposition} 

{\em Proof.} 
We need some material concerning amalgamation bases in 2-step nilpotent groups, see Theorems 3.3 - 3.6 in \cite{sar}. 
It can be deduced from these theorems  that a 2-step nilpotent group $G$ of exponent $p$ is an amalgamation base if and only if $G'= Z(G)$. 
Below we only need sufficiency of this for finite groups, and the latter can be shown easier. 
Indeed, consider a finite $G \in \mathbb{G}_{2_{nil}p}$ with $G'=Z(G)$. 
Note that any embedding of $G$ into any $H \in \mathbb{G}_{2_{nil}p}$ preserves every unary predicate $P$ with the property $H' \le P \le Z(H)$ where $P$ on $G$ is defined by $P=G' = Z(G)$. 
By Lemmas 2.2 and 2.3 of \cite{baud} the class 
of finite members of $\mathbb{G}_{2_{nil}p}$ in the language extended by $P$ with $H' \le P \le Z(H)$, has the amalgamation property. 

It now suffices to verify that every finite group $G_0 \in \mathbb{G}_{2_{nil}p}$ extends to a finite group $G_1 \in \mathbb{G}_{2_{nil}p}$ with $G'_1 = Z(G_1 )$. 
This can be again deduced from the amalgamation result of \cite{baud}.  
Indeed, consider $G_0$ with the predicate $P(G_0 )=Z(G_0)$. 
Since $P(G_0)$ is a vector space over $\mathbb{F}_p$ we decompose $P(G_0) = G'_0 \oplus \langle e_1 ,\ldots , e_k \rangle$ where $e_1 ,\ldots , e_k$ is a basis of the complement of $G'_0$. 
Now take a finite $H\in \mathbb{G}_{2_{nil}p}$ such that 
$H' = Z(H) = \langle e_1 ,\ldots , e_k \rangle$. 
The amalgamation process described in Lemmas 2.2 and 2.3 of \cite{baud} applied to $(G_0,P)$ and $(H,P)$ gives the result. 
$\Box$ 

\bigskip 

We do not know if the statement of this proposition holds for $c$-step nilpotent groups where $c>2$. 
The argument in the beginning of the proof together with \cite{baud04} proves that if 
a finite group $G\in \mathbb{G}_{c_{nil}p}$ satisfies  $\Gamma_i (G) = Z_{c+1 - i} (G)$ for all $1 < i \le c$, then  $G$ is an amalgamation base. 
{\em Is it possible to embed every finite $H\in \mathbb{G}_{c_{nil}p}$ into some $G$ with this property?}

\begin{remark} 
{\em 
The question formulated in the previous line has a strong connection with the following issue. Corollaries 1.3 - 1.5 of \cite{baud04} state that the class of all finite $L^{+}$-structures from  $\mathbb{G}^P_{c,p}$ has a Fra\"{i}ss\'{e} limit, say $\mathsf{G}^+_{c,p}$, and it is $\aleph_0$-categorical, ultrahomogeneous and universal for countable groups from $\mathbb{G}_{c_{nil}p}$. 
It is natural to ask whether the reduct of this structure, say $\mathsf{G}_{c,p}$, to the group language is a generic structure from Theorem \ref{cnilp}. 
The following sketch gives the affirmative under the assumption that the question formulated before the remark has the positive answer. 
According to Lemma 1.1 from \cite{iva99} it suffices to prove that the reduct  $\mathsf{G}_{c,p}$ belongs to every set of the form {\bf S}, {\bf J}, {\bf E} from that lemma. 
The cases of {\bf J} and {\bf E} are rather easy by residual finiteness of finitely generated nilpotent groups and ultrahomogeneity of  $\mathsf{G}^+_{c,p}$.  
We remind the reader that given a diagram $D$ the set ${\bf S}(D)$ consists of all enumerated groups $\mathsf{G}\in \mathcal{G}_{c_{nil}p}$ such that when $\mathsf{G}$ realizes $D$ then it also realizes some extension of $D$ witnessing WAP over $D$ in $\mathcal{B}_{c_{nil}p}$.  
To see this case it suffices to show that any finite subgroup $H_1 <\mathsf{G}_{c,p}$
extends to a finite $H_2 < \mathsf{G}_{c,p}$
such that every pair of finite extentions of $H_2$ is amalgamated over $H_1$.  
Take a finite $G\in \mathbb{G}_{c_{nil}p}$ that extends $H_1$ and satisfies  $\Gamma_i (G) = Z_{c+1 - i} (G)$ for all $1 < i \le c$. 
By ultrahomogeneity of $\mathsf{G}^+_{c,p}$, the unique expansion $G^+$ of $G$ to a structure from $\mathbb{G}^P_{c,p}$ has a copy in  $\mathsf{G}^+_{c,p}$ over $H_1$. 
Let $H_2$ be the reduct of this copy to the group language. 
Since it is an amalgamation base, we can finish the argument.  
This, in particular, shows that the group $\mathsf{G}_{2,p}$ is a generic structure of the space $\mathcal{G}_{2_{nil}p}$ (see Proposition \ref{nil-cap}). 
} 
\end{remark}

\section{Strong undecidability implies non-WAP} 

In this section we give several group-theoretic properties such that $\mathcal{X}_{M,T}$ and the corresponding space of enumerated groups does not have generics. 
Thus, the section will focus on cases {\em I} - {\em II} of Section 1.2. 
On the other hand, the main technical tool, Theorem \ref{nWAP}, works in much more general situations.

\subsection{Completeness and Kuznetsov's theorem} 

In Section 3 we assume that $M$ is an $\omega$-categorical  structure in a language $L$ with a distinguished sort $\mathsf{Gp}$, and  $T$ is an extension of $Th(M)$ in some language $L\cup \bar{r}$ satisfying the assumptions of Section 1.1 
(i.e. our assumptions on $M$ are slightly stronger than in Section 2.1).  
When $(M,{\bf \bar{r}})$ is an expansion of $M$ in this language, we consider $\mathsf{Gp}$ as a structure with respect to those relations/functions of $L\cup \bar{r}$ which do not depend on elements of other sorts. 
Then any subset of $\mathsf{Gp}$ closed under 
operations of $\mathsf{Gp}$ is viewed as a substructure. 
If such a substructure is generated by some finite set $A$, then all valid atomic formulas depending on tuples from $A$ form a presentation of this substructure.  
Further notation corresponds to Section 2.1.

The following proposition can be viewed as a generalization of the theorem of Kuznetsov that a recursively presented simple group has decidable word problem. 

\begin{proposition} \label{Kuzn}  
Assume that $D(\bar{c}) \in \mathcal{B}_T$ and  an extension $D'(\bar{c}') \in \mathcal{B}_T$ is t-isolating for $\bar{c}$. 
Let $(M, {\bf \bar{r}})\in [D'(\bar{c}')]$. 
Assume that $Th(M, \bar{c}') \cup D'(\bar{c}')\cup T$ is a computably enumerable set. 
Then in $(M, {\bf \bar{r}})$, the elements of $\bar{c}$ that belong to  the sort $\mathsf{Gp}$ generate a recursively presented structure with decidable theory $\Xi_{\, \bar{c}\cap \mathsf{Gp}, \bf \bar{r}}$.  
\end{proposition} 

{\em Proof.} 
We may assume that $\bar{c}$ contains all elements of the sort $\mathsf{Gp}$ which are distinguished in $(M, {\bf \bar{r}})$ as constants.  
Let us extend $D(\bar{c})$ by all atomic formulas with parameters just from $\bar{c}\cap \mathsf{Gp}$ that can be logically deduced from $Th(M, \bar{c}') \cup D'(\bar{c}')\cup T$.  
We also include negations of such formulas when they are deduced from $Th(M, \bar{c}') \cup D'(\bar{c}')\cup T$. 
The obtained set, say $\Sigma$, is computably enumerable. 
The set of atomic formulas from $\Sigma$ forms a computably enumerable presentation in the generators $\bar{c} \cap \mathsf{Gp}$ of the corresponding substructure generated in $\mathsf{Gp}$.
Remember, that by t-isolation of $D'(\bar{c}')$ for $\bar{c}$, this presentation does not depend on the expansion $(M, {\bf \bar{r}})$.  

The deciding algorithm is as follows. 
Let $\phi$ be an atomic formula with parameters from $\bar{c} \cap \mathsf{Gp}$.  
We verify $\phi \in \Sigma$ using the enumeration of $\Sigma$ obtained above.  
The algorithm stops when this procedure gives $\phi$ or $\neg \phi$ in $\Sigma$. 
If the first case does not hold, then in any member of $[D'(\bar{c}')]$ the substructure generated by $\bar{c} \cap \mathsf{Gp}$ satisfies $\neg \phi$. 
If $\neg \phi$ does not belong to $\Sigma$, then $Th(M, \bar{c}') \cup D'(\bar{c}') \cup T\cup \{ \neg \phi \}$ is consistent. 
Any model of this theory does not have a reduct that is isomorphic to $(M, \bar{c}')$ (since otherwise it defines a member of $[D'(\bar{c}')]$). 
This contradicts $\omega$-categoricity of $Th(M,\bar{c}')$. 
We see that the second case holds. 
$\Box$ 

\bigskip  

The following theorem will be the main tool for proving the absence of generics in $\mathcal{G}_T$.

\begin{theorem}  \label{nWAP} 
Assume that $D(\bar{c})\in \mathcal{B}_T$ (where $\bar{c}$ includes all distinguished elements of $M$) satisfies the property that every extension $D'(\bar{c}')$ defines with $T \cup Th(M,\bar{c}')$ a computably enumerable theory, and for every $(M, {\bf \bar{r}})\in [D(\bar{c})]$ the substructure of $\mathsf{Gp}$ generated by $\bar{c} \cap \mathsf{Gp}$ in $(M, {\bf \bar{r}})$ has undecidable theory $\Xi_{\, \bar{c}\cap \mathsf{Gp}, \bf \bar{r}}$. 

Then $D(\bar{c})$ does not have an extension required by the property that $\mathcal{B}_T$ satisfies WAP. 
\end{theorem} 

{\em Proof. } 
Let $D'(\bar{c}')$ be an extension of $D(\bar{c})$. Let us show that it does not witness WAP. 
By Proposition \ref{Kuzn}, $D'(\bar{c}')$ is not t-isolating for $\bar{c}$. 
This means that there are $(M, {\bf \bar{r}})$, $(M, {\bf \bar{r}'}) \in [D'(\bar{c}')]$ and an atomic formula $\phi (\bar{c}_0)$ with 
$\bar{c}_0 \subseteq \bar{c} \cap \mathsf{Gp}$ 
such that 
$(M, {\bf \bar{r}}) \models \phi(\bar{c}_0)$ but 
$(M, {\bf \bar{r}'}) \models \neg \phi(\bar{c}_0)$. 
We see that the two extensions of $D'(\bar{c}')$ by $\phi (\bar{c}_0)$ and by $\neg \phi(\bar{c}_0)$ (where domains of these extensions are enriched by values of all terms appeared in $\phi(\bar{c}_0)$) cannot be amalgamated over $D(\bar{c})$. 
Therefore, $D'(\bar{c}')$ does not witness WAP. 
$\Box$ 

\bigskip 

In the following corollary we assume Scenarios
{\em I - II} of Section 1.2. 

\begin{corollary} \label{nWAPgr} 
Assume that $D(\bar{c})\in \mathcal{B}_T$ (where $\bar{c}$ includes all distinguished elements of $M$) satisfies the property that every extension $D'(\bar{c}')$ defines with $T \cup Th(M,\bar{c}')$ a computably enumerable theory, and for every $(M, {\bf \bar{r}})\in [D(\bar{c})]$ the subgroup generated by $\bar{c} \cap \mathsf{Gp}$ in $(M, {\bf \bar{r}})$ has undecidable word problem. 

Then $D(\bar{c})$ does not have an extension required by the property that $\mathcal{B}_T$ satisfies WAP. 
\end{corollary}

\begin{remark}
{\em 
Corollary \ref{nWAPgr} can be viewed as a counterpart of Theorem 1.1.6 of \cite{GEL}. 
Indeed, by Theorem 1.2 from \cite{iva99} both statements have the same conclusion: the absence of generic groups in $\mathcal{G}_T$. 
Their assumptions are also similar: in each of them undecidability of the word problem is involved in some hereditary way. 
As we will see below, both statements can be applied to the same objects with the same result. 
We also mention that Theorem 1.1.6 of \cite{GEL} is a result of substantial effort: it appears in the final section of this long paper and is based on Theorem 1.1.3, another main theorem of the paper
(among four main theorems which are stated by the authors in introduction). 
} 
\end{remark}

We now give several applications of Corollary \ref{nWAPgr} for concrete spaces of enumerated groups. 
In these applications we are under Scenario {\em I} of Section 1.2. 
In particular, the structure $M$ is just $(\omega,1)$. 

\begin{corollary} \label{darbi} 
The space of torsion free enumerated groups $\mathcal{G}_{t.f}$ does not have a generic group. 
\end{corollary} 

{\em Proof.} 
Since the class of torsion free groups is universally axiomatizable, all torsion free groups form a closed subspace of the logic space $\mathcal{G}$ of enumerated groups. 
Furthermore, this axiomatization is computable. 
In particular every member of $\mathcal{B}_{T}$ is computably enumerable. 

In Proposition 5.2.13 of \cite{GEL} the authors present a construction (of A. Darbinyan) of a finitely presented left orderable group $G = \langle \bar{c} \, | \, \mathcal{R}(\bar{c}) \rangle$ with two elements $w_1 (\bar{c}), w_2 (\bar{c})\in G$ such that every homomorphic image of $G$ where $w_1 (\bar{c})\not= w_2 (\bar{c})$, has undecidable word problem. 
This group is torsion free.  
We may assume that values of all subwords appearing in $\mathcal{R}(\bar{c}) \cup \{ w_1 (\bar{c}), w_2 (\bar{c})\}$ are in $\bar{c}$. 
Applying Corollary \ref{nWAPgr} we see that $(\bar{c}, \mathcal{R}(\bar{c}) \cup \{ w_1 (\bar{c}) \not= w_2 (\bar{c}) \} )$ is a diagram without extensions witnessing WAP. 
Indeed, for every $(M, {\bf \bar{r}})\in [D(\bar{c})]$  the group $\langle \bar{c} \rangle$ (defined in $(M, {\bf \bar{r}})$) is a homomorphic image of $G$ and therefore has an undecidable word problem. 
Applying Theorem 1.2 from \cite{iva99} we get the result. 
$\Box$ 

\bigskip 

Note that the proof of this corollary also demonstrates that the space $\mathcal{G}$ of all enumerable groups does not have generic groups. 
It also proves that the subspace of left orderable groups does not have generics too. 
These statements are among the major theorems of \cite{GEL}. 
On the other hand, we should emphasize that the idea of Darbinyan's group used above is taken from \cite{GEL}.  
We now give a new application. 
It is in some way opposite to the material in Section 2.3.  
Let $\mathbb{G}_{exp.n}$ denote the variety of groups of exponent $n$. 

\begin{corollary} \label{burnside} 
There is a constant $C$ such that for every odd integer $n \ge C$, the space $\mathcal{G}_{exp.n}$ of all groups of exponent $n$ does not have a generic group. 
\end{corollary} 

{\em Proof.} 
Since the class of torsion groups of exponent $n$ is universally axiomatizable, this class forms a closed subspace of the logic space of enumerated groups. 
Furthermore, this axiomatization is computable. 
In particular, every member of $\mathcal{B}_{T}$ is computably enumerable. 

It is proved in \cite{OlshanskiiJA} that there is a constant $C$ such that for every odd integer $n \ge C$,  the following
is true.
Let a group $H$ from the variety $\mathbb{G}_{exp.n}$ be given by a countable set of generators $a_1, a_2, \ldots$ 
and a recursively enumerable set of defining relations $\mathcal{R}$ in these generators. 
Then $H$ is a subgroup of a 2-generated group $E$ finitely presented in $\mathbb{G}_{exp.n}$.
Furthermore, it easily follows from the description given in \cite{OlshanskiiJA} that the embedding $H \rightarrow E$ is computable. 
Indeed, it is shown in the proof of Corollary 1.2 of \cite{OlshanskiiJA} that the group $H$ 
can be presented as a subgroup of $E=B(2,n)/P$ 
of the form $KP/P$, where 
\begin{itemize} 
\item $B(2,n)$ is the free 2-generated Burnside group, 
\item $K \cong B(\infty , n)$  the free $\omega$-generated Burnside group, 
\item $K = \langle w_1, w_2, \ldots \rangle$ and  $w_1, w_2 , \ldots$ is a computable set of words in $B(2,n)$,    
\item the set $P$ in this presentation is computably enumerable. 
\end{itemize} 
Thus, the map $a_i \rightarrow w_i$ defines the corresponding embedding in a computable way. 

We now perform the Miller's construction (i.e. from Proposition 5.2.13 of \cite{GEL}) as follows.  
Fix recursively inseparable recursively enumerable subsets $N,N' \subset \omega$ (i.e. there is no recursive set $W$ such that $N \subset W$ and $N'\cap W = \emptyset$ ). 
Let $A_0$ be the free abelian group on countably many generators $a_i$. 
Consider the following abelian quotient $A$ of $A_0$: 
\[ 
A = 
\langle a_1 , a_2 , \ldots \, | \, 
a^n_i = 1 \, (i \in \omega )\, , \, 
a_1 = a_i \mbox{ if } i \in N \mbox{ and } 
a_2 = a_j \mbox{ if } j \in N' \rangle.
\]  
Let $\pi_0$ be the quotient map from
$A_0$ to $A$. 
Applying the construction of Corollary 1.2 of \cite{OlshanskiiJA} (described in the beginning of our proof) we computably embed $A$ into a 
2-generated group $G$, which is finitely presented in $\mathbb{G}_{exp.n}$.  
Denote this embedding by $\phi$.  
Set 
$F = \{ \phi( a_1 ), \phi (a_2)\} \subset G$.  Let us verify that, for any surjective homomorphism, say $\pi$ from $G$ to a group
$H$ that is injective when restricted to
$F$, we have that $H$ has undecidable word problem. 
Indeed, if the word problem of $H$
was decidable, then the set
$W := \{ i \in \omega \, | \,  
(\pi \phi \pi_0)(a_i)=(\pi \phi \pi_0 )(a_1) \}$ 
is computable, contains $N$, and is disjoint from $N'$, contradicting the recursive inseparability of $N$ and $N'$. 

We may assume that $G$ is presented in $\mathbb{G}_{exp.n}$ as $\langle \bar{c}\, |\, \mathcal{R}(\bar{c}) \rangle$ where values of all subwords appearing in the presentation are in $\bar{c}$. 
Thus $\{ v(\bar{c}') = 1 \, | \, v(\bar{c}') \in \mathcal{R}(\bar{c}) \} \cup \{ (\phi (a_1 ) \not= \phi (a_2 ) \}\in \mathcal{B}_T$ 
where $T$ is the theory of all periodic grous of exponent $n$,  
and this diagram does not have extensions witnessing WAP. 
The corresponding argument is as in the final part of the proof of Corollary \ref{darbi}.   
$\Box$ 

\bigskip 

Comparing this result with Section 2.3 we see how crucial the property of nilpotence is in the context of enumerated groups.

\subsection{Generics over rationals} 

We now give an application of Theorem \ref{nWAP} to actions on $\mathbb{Q}$ (and thus on the real line $\mathbb{R}$). 
Consider Scenario {\em II} of Section 1.2, where $M_0 = (\mathbb{Q}, < )$. 
We remind the reader that this structure is $\omega$-categorical, admits elimination of quantifiers and its theory is decidable.  
Furthermore, these properties still hold for any expansion $(\mathbb{Q},< , \bar{a})$ by a tuple $\bar{a} \in \mathbb{Q}$. 
Note that the  theory of this expansion is determined by the size of $\bar{a}$.   

Now define $M$ to be $M_0$ with the additional sort $\mathsf{Gp}$ and the constant symbol $1$ interpreted by $1\in \omega$. 
The symbols $\bar{r}$ are $\cdot$,  $^{-1}$ and $\mathsf{ac}$, for an action $M_0\times \mathsf{Gp}\rightarrow M_0$ by automorphisms. 
 
The theory $T$ contains the universal axioms of groups on the sort $\mathsf{Gp}$ with the unit $1$ and the axioms for an action by automorphisms  (i.e. $<$ is preserved).   
The space $\mathcal{X}_{M,T}$ is the logic space of all countable expansions of $M$ where the group structure is defined on $\mathsf{Gp}$, with a fixed unit $1$ and an action $\mathsf{ac}$ . 
We denote it by $\mathcal{X}_{(\mathbb{Q},<)}$. 

\begin{corollary} \label{Q} 
The space $\mathcal{X}_{(\mathbb{Q},<)}$ does not have a generic action. 
\end{corollary} 

{\em Proof.} 
It is easy to see that $\mathcal{X}_{(\mathbb{Q},<)}$ is a closed subspace of the logic space of group expansions of $M$ as above. 

As in Corollary \ref{darbi} we use the  construction of A. Darbinyan of a finitely presented left orderable group $G = \langle \bar{c} | \mathcal{R}(\bar{c}) \rangle$ with two elements $w_1 (\bar{c}), w_2 (\bar{c})\in G$ such that any homomorphic image of $G$ where $w_1 (\bar{c})\not= w_2 (\bar{c})$, has undecidable word problem  
(see Proposition 5.2.13 of \cite{GEL}). 
 
The group $G$ acts on $(\mathbb{Q},<)$ by automorphisms.  
Indeed, take the ordering $<$ on $G$ so that $G$ acts on it from the left. 
If $<$ is dense we identify it with $\mathbb{Q}$. 
If $<$ is discrete then embed it into $(\mathbb{Q},<)$. 
Its convex hull in $\mathbb{Q}$ can be identified with $(\mathbb{Q}, <)$ again (it is easy to see that there is no maximal element in $(G,<)$).  
As a result we decompose $(\mathbb{Q}, <)$ into a family of intervals with endpoints in $G$.  
All of them can be identified with the interval $(1, g_0)$, where $g_0$ is the least element of $G$ greater than $1$. 
Extend the action of $G$ on $G$ (by left multiplication) to the rest of $\mathbb{Q}$
using this identification. 
   
Now by Corollary \ref{nWAPgr} we see that $\mathcal{R}(\bar{c}) \cup \{ w_1 (\bar{c}) \not= w_2 (\bar{c}) \}$ defines a diagram without extensions witnessing WAP. 
Applying Theorem 1.2 from \cite{iva99} we get the result. 
$\Box$

\bigskip 

This result fits to the topic of generic representations of countable groups in Polish groups, for example see \cite{DM}. 
A typical situation studied in these investigations is as follows. 
Take a countable group $\Gamma$ and a countable first order structure $M$. 
Now let $\mathsf{Rep} (\Gamma ,\mathsf{Aut} (M))$  be the space of all representations of $\Gamma$ in  $\mathsf{Aut} (M)$. 
It is considered under the natural Polish topology (from $\mathsf{Aut}(M)^{\Gamma}$) and the conjugacy action of $\mathsf{Aut} (M)$. 
The notion of a generic representation in $\mathsf{Rep} (\Gamma ,\mathsf{Aut} (M))$ is defined in the standard way, see Introduction in \cite{DM}.  
Note that every representation of $\Gamma$ can be viewed as an expansion $(M, g_{\gamma})_{\gamma \in \Gamma}$ where each $g_{\gamma}$ is the  automorphism corresponding to $\gamma$.  
Thus, under Scenario {\em II} of Section 1.2, $\mathsf{Rep} (\Gamma ,\mathsf{Aut} (M))$ forms the subset of $\mathcal{X}_{M,T}$ where the group structure of $\mathsf{Gp}$ is fixed (and coincides with $\Gamma$). 
It is well known that $\mathsf{Rep} (\mathbb{Z} ,\mathsf{Aut} (\mathbb{Q},<))$ has a generic representation (see \cite{truss}), but $\mathsf{Rep} (\mathbb{F}_2 ,\mathsf{Aut} (\mathbb{Q},<))$ does not have it (an unpublished theorem of I.Hodkinson).  
Corollary \ref{Q} naturally complements these facts. 
Furthermore, it seems plausible that the fact concerning $\mathbb{F}_2$ can be useful for an alternative proof of it. 

\section{Semigroups and rings} 

The approach of Section 3 can be applied to other structures studied in algebra. 
In this section we consider how it works in the cases of semigroups and rings. 

In the cases of semigroups we modify Scenario {\em I} of Section 1.2 so that the structure $M$ is just $\omega$ and $L$ is empty (i.e. in terms of Section 1.2 the sort $\mathsf{Gp}$ coincides with $M$).  
Let $\bar{r}$ consist of the binary function of multiplication \, $\cdot$. 
Let $T$ be the universal theory of semigroups.  
Under this convention each element of $\mathcal{B}_T$ consists of a tuple $\bar{c}\subset \omega$ and a description of a partial function for $\cdot$ on $\bar{c}$ (the associativity  axiom for $\cdot$ is also added), i.e. a partial semigroup. 
The space $\mathcal{X}_{M,T}$ is the logic space of all countable semigroups. 
We denote it by $\mathcal{S}$. 

The case of rings is similar. 
We now assume that $\bar{r}$ consist of three binary functions of addition \, $+$ \, , subtraction \, $-$ \, and multiplication \, $\cdot$ . 
We also include the constant symbol $0$ into $\bar{r}$. 
Let $T$ be the universal theory of associative rings (without unit).  
Thus, each element of $\mathcal{B}_T$ consists of a tuple $\bar{c}\subset \omega$ and a description of three  partial functions for $+$, $-$ and $\cdot$ on $\bar{c}$ (with $0$ and axioms of a ring). 
The space $\mathcal{X}_{M,T}$ corresponding to $T$ will be denoted by $\mathsf{Ring}$. 

\bigskip 

{\em Completeness and Kuznetsov's theorem for semigroups and rings} work as follows. 
It is obvious that in both cases Definition \ref{isol} is equivalent to the following one. 

\begin{definition}
Let $\bar{c}$ be a subtuple of $\bar{c}_1$ and $D(\bar{c}_1) \in \mathcal{B}_T$. 
We say that $D(\bar{c}_1)$ is {\em t-isolating} for $\bar{c}$ if in every member of $[D(\bar{c}_1)]$ the subsemigroup (resp. subring) generated by $\bar{c}$ in $M$ satisfies the same semigroup (resp. ring) equalities of the form $w_1 (\bar{c}') = w_2 (\bar{c}'')$ with $\bar{c}', \bar{c}'' \subseteq \bar{c}$. 
\end{definition}
  
%
%
%
The results of Section 3 have natural counterparts for semigroups and rings. 
For example, the following statement is a counterpart of Corollary \ref{nWAPgr}. 

\begin{corollary}  \label{nWAPsemi} 
Assume that $D(\bar{c})\in \mathcal{B}_T$ satisfies the property that for every $(M, {\bf \cdot})$ (resp. $(M, + , 0 , - , \cdot )$) $\in [D(\bar{c})]$ the subsemigroup (resp. subring) generated by $\bar{c}$ in $(M, {\bf \cdot})$ (resp. $(M, + , 0 , - , \cdot )$) has undecidable word problem. 

Then $D(\bar{c})$ does not have an extension required by the property that $\mathcal{B}_T$ satisfies WAP. 
\end{corollary} 

The easy form of this corollary follows from the fact that $T$ is computably enumerable in it, and the language of $M$ is empty.  
Corollary \ref{nWAPsemi} is applied to  enumerated semigroups as follows.

\begin{corollary} \label{darbisemi} 
The space $\mathcal{S}$ of enumerated semigroups does not have a generic semigroup. 
\end{corollary} 

{\em Proof.} 
We repeat the argument of Corollary \ref{darbi}. 
Let $G = \langle \bar{c} \, | \, \mathcal{R}(\bar{c}) \rangle$ be a finitely presented left orderable group  with two elements $w_1 (\bar{c}), w_2 (\bar{c})\in G$ such that any homomorphic image of $G$ where $w_1 (\bar{c})\not= w_2 (\bar{c})$, has undecidable word problem  
(see Proposition 5.2.13 of \cite{GEL}). 
This group can be converted into a finitely presented semigroup as follows.  
Let $n= |\bar{c}|$. 
For every $c_i \in \bar{c}$ we introduce the letter $c_{n+i}$. 
The group $G$ has a semigroup presentation for generators 
$\bar{c}, 1, c_{n+1},\ldots , c_{2n}$ as follows. 
We start with relations saying that $1$ is the unit: $c_i \cdot 1 =  1 \cdot c_i = c_i$, $1 \le i \le 2n$ and $1\cdot 1 = 1$. 
Now, add all relations of the form $c_i \cdot c_{n+i} = c_{n+i} \cdot c_i =1$, $1 \le i \le n$. 
Every word from $\mathcal{R}(\bar{c})$ 
can be rewritten as a semigroup word by replacing each $c^{-1}_i$ by $c_{n+i}$. 
We extend the family of relations by all equalities 
$w(\ldots, c_i,\ldots c_{n+i}, \ldots )=1$ for all words rewriting ones from $\mathcal{R}(\bar{c})$.  
It is easy to see that the semigroup with this presentation is isomorphic to $G$. 
Let $w'_1$ and $w'_2$ be obtained from 
$w_1 (\bar{c})$ and $w_2 (\bar{c})$ using the rewriting procedure described above for $\mathcal{R}(\bar{c})$. 
Applying Corollary \ref{nWAPsemi} we see that the semigroup presentation of $G$ together with $w'_1  \not= w'_2$ is a diagram without extensions witnessing WAP. 
This follows from the argument of Corollary \ref{darbi} together with the  observation that every semigroup homomorphism from $G$ is a group homomorphism. 
The rest is clear. 
$\Box$ 

\bigskip 

Let us consider the case of associative rings. 

\begin{corollary} \label{ring} 
The space $\mathsf{Ring}$ of all assoiative rings does not have a generic ring. 
\end{corollary} 

{\em Proof.} 
V. Belyaev has proved in \cite{bel} that given a countable associative ring $A$ and an enumeration $A= \{ a_1 , \ldots , a_i , \ldots \}$ (where repetitions are possible) such that $A$ is presented by a computably enumerable set of relations in this alphabet, there is a finitely presented ring $R$, say with generators $\{ a, b, c, d_1 , \ldots , d_k\}$, such that the function  
$a_i \rightarrow ab^i c$, $1 \le i$,  embeds $A$ into $R$. 
 
We perform the Miller's construction (i.e. from Proposition 5.2.13 of \cite{GEL}) in the case of rings as follows.  
Fix recursively inseparable recursively enumerable subsets $M,N \subset \omega$. 
Let $A_0$ be the ring $\oplus^{\omega}_{i=1} \mathbb{Z}$ of countably many free generators $a_i$, $i\in \omega$. 
Consider the following quotient $A$ of $A_0$: 
\[ 
A = 
\langle a_1 , a_2 , \ldots \, | \,  
a_1 = a_i \mbox{ if } i \in M \mbox{ and } 
a_2 = a_j \mbox{ if } j \in N \rangle.
\]  
Let $\pi_0$ be the quotient map from
$A_0$ to $A$. 
Applying the construction of Belyaev from \cite{bel} (described in the beginning of our proof) we computably embed $A$ into a finitely generated ring $R$, which is finitely presented in $\mathsf{Ring}$.  
Denote this embedding by $\phi$.  
Set 
$F = \{ \phi( a_1 ), \phi (a_2)\} \subset R$.  
Repeating the argument of Corollary \ref{burnside} it is easy to verify that, for any surjective homomorphism, say $\pi$ from $R$ to a ring 
$R'$ that is injective when restricted to $F$, we have that $R'$ has undecidable problem of equality of terms. 

We see that extending the presentation of $R$ by $\phi (a_1 ) \not= \phi (a_2 )$ (assuming that the values of subterms of terms of the presentation  are among generators) we obtain a diagram without extensions witnessing WAP. 
The corresponding argument is as in the final part of the proof of Corollary \ref{darbi}.   
$\Box$

\vspace*{10mm}

\begin{flushleft}
\begin{footnotesize}
Aleksander Iwanow, \, Inst. Comp. Sci., University of Opole, ul. Oleska 48, 45-052 Opole, Poland\\
email: aleksander.iwanow@uni.opole.pl
\end{footnotesize} 

\bigskip 

\begin{footnotesize}
Krzysztof Majcher, \, Faculty of Information and Communication Technologye,
Wroc{\l}aw University of Technology, \,
pl. Grunwaldzki 13, \, 
50--377 Wroc{\l}aw, 
Poland\\
email: k.majcher@pwr.edu.pl
\end{footnotesize}

\end{flushleft}


\begin{thebibliography}{20} 

\bibitem{baud} A. Baudisch, A Fra\"{i}ss\'{e} limit of nilpotent groups of finite exponent, Bull. London Math. Soc. 33 (2001) 513 -- 519. 
\bibitem{baud04} A. Baudisch, More Fra\"{i}ss\'{e} limits of nilpotent groups of finite exponent, Bull. London Math. Soc. 36 (2004) 613 -- 622. 
\bibitem{bel} V.Ya. Belyaev, Subrings of finitely presented associative rings. 
Algebra and Logic 17 (1978) 407 -- 414 (1979) (Russian version: 
Algebra i Logika 17 (1978), no. 6, 627 -- 638, 746). 
\bibitem{CMZ} A. E. Clement, S. Majewicz and M. Zyman, The Theory of Nilpotent Groups, Springer International Publishing AG 2017. 
\bibitem{DM} M. Doucha and M. Malicki, Generic representations of countable groups, Trans. Amer. Math. Soc. 372 (2019) 
8249 -- 8277. 
\bibitem{marton}  M. Elekes, B. Geh\'{e}r, K. Kanalas, T. K\'{a}tay and T. Keleti, Generic countably infinite groups, 
arXiv:2110.15902, 2021. 
\bibitem{marton-top} M. Elekes, B. Geh{\'e}r, T. K{\'a}tay, T. Keleti, A. Kocsis and M. P{\'a}lfy, Generic properties of topological groups,
Proc. Royal Soc. Edinburgh, Section A: Mathematics, 2024, DOI:https://doi.org/10.1017/prm.2024.91. 
\bibitem{GEL} 
I. Goldbring, S. Kunnawalkam Elayavalli, Y. Lodha, Generic algebraic properties in spaces of enumerated groups, Trans. Amer. Math. Soc.	
376 (2023) 6245 -- 6282. 	
\bibitem{iva99} A. Ivanov, Generic expansions of $\omega$-categorical structures and semantics of generalized quantifiers, J. Symb. Logic 64 (1999) 775 -- 789. 
\bibitem{KR} A. S. Kechris and Ch. Rosendal, Turbulence amalgamation and generic automorphisms of homogenoeus structures, Proc. Lond. Math. Soc. (3) 94 (2007) 302 -- 350.  
\bibitem{keisler} H.J. Keisler, Fundamentals of Model Theory,   
in: Handbook of Mathematical Logic (J. Barwise ed.) pp. 47 -- 103, North-Holland, NY, 1977. 	
\bibitem{KKKP} A. Krawczyk, A. Kruckman, W. Kubi\'{s} and A. Panagiotopoulos, Examples of weak amalgamation classes, Math. Logic Quart. 68(2) (2022) 178 -- 188.	
\bibitem{MI} B. Majcher-Iwanow, $G_{\delta}$-pieces of canonical partitions of $G$-spaces, Math Logic. Quart. 51 (2005) 450 -- 561. 
\bibitem{OlshanskiiJA} A. Yu. Olshanskii, Groups finitely presented in Burnside varieties, J. Algebra 560 (2020) 960 -- 1052. 
\bibitem{sar} D. Saracino, 
Amalgamation bases for nil-2 groups, Alg. Univ. 16 (1983) 47 -- 62.  
\bibitem{truss} J. Truss, Generic automorphisms of homogeneous structures, Proc. London Math. Soc. (3) 65(1992) 121--141.
\end{thebibliography}
\end{document}